\documentclass[11pt, oneside]{article}   	% use "amsart" instead of "article" for AMSLaTeX format
\usepackage{geometry}                		% See geometry.pdf to learn the layout options. There are lots.
\usepackage{graphicx}				% Use pdf, png, jpg, or eps with pdflatex; use eps in DVI mode
\usepackage{latexsym}
\usepackage{amsmath}
\usepackage{amsthm}
\usepackage{amscd}
\usepackage{amssymb}
\usepackage{mathrsfs}
\usepackage{mathabx}
\usepackage{txfonts}

\newtheorem{theorem}{Theorem}[subsection]

\theoremstyle{definition}
\newtheorem{definition}[theorem]{Definition}

\theoremstyle{remark}

\newcommand{\R}{\mathrm{R}}

\newcommand{\tr}{\mathrm{tr}}

\renewcommand{\d}{\mathrm{d}}

\newcommand{\slashepsilon}{\epsilon\llap\slash}

\newcommand{\circDelta}{\stackrel{\circ}{\Delta}}
\numberwithin{equation}{section}
\numberwithin{theorem}{section}

\DeclareMathAlphabet{\mathpzc}{OT1}{pac}{m}{it}

\title{The Intersection of a Hyperplane with a Lightcone in the Minkowski Spacetime}
\author{Pengyu Le}
\date{}

\begin{document}
\maketitle

\begin{abstract}
Klainerman, Luk and Rodnianski derived an anisotropic criterion for formation of trapped surfaces in vacuum, extending the original trapped surface formation theorem of Christodoulou. The effort to understand their result led us to study the intersection of a hyperplane with a lightcone in the Minkowski spacetime. For the intrinsic geometry of the intersection, depending on the hyperplane being spacelike, null or timelike, it has the constant positive, zero or negative Gaussian curvature. For the extrinsic geometry of the intersection, we find that it is a noncompact marginal trapped surface when the hyperplane is null. In this case, we find a geometric interpretation of the Green's function of the Laplacian on the standard sphere. In the end, we contribute a clearer understanding of the anisotropic criterion for formation of trapped surfaces in vacuum.
\end{abstract}

\section{Introduction}

In 1965 work of Penrose \cite{P}, he introduced the concept of a closed trapped surface, i.e. a closed spacelike surface with the property that the area element decreases pointwise  for any infinitesimal displacement along the future null direction normal to the surface. Based on this concept, Penrose proved his famous incompleteness theorem. The incompleteness theorem says the following: a vacuum spacetime cannot be future null geodesically complete if it contains a non-compact Cauchy hypersurface and moreover contains a closed trapped surface. Then it follows that the spacetime with a complete null infinity which contains a closed trapped surface must contain a future event horizon the interior of which contains the closed trapped surface.

The incompleteness theorem gives rise to the following problem: could the closed trapped surface be formed dynamically in vacuum? In the important monograph of Christodoulou \cite{C},  he answered this problem positively by showing that the closed trapped surface could be formed by sufficiently strong incoming gravitational waves. Christodoulou introduced a new method called ``short pulse method", which allows him to establish a longtime existence theorem for a development of initial data which is large enough so that the trapped surfaces have a chance to form within this development.

Klainerman, Luk and Rodnianski derived an anisotropic criterion for formation of trapped surfaces in vacuum in \cite{K-L-R}. The original trapped surface formation theorem requires a lower bound on energy per solid angle of the incoming gravitational waves in all directions. They showed instead that, if the energy per solid angle of the incoming gravitational waves in some neighbourhood of some direction is sufficiently large depending on the size of the neighbourhood, then the development of initial data contains a trapped surface. Their result is surprising since the spacetime could be flat in a large part of the development of initial data they considered.

One of the key ideas in \cite{K-L-R} is that given the existence result of \cite{C}, one needs to solve an elliptic inequality to find a surface in the spacetime which is trapped.

We find that the construction of the trapped surface in \cite{K-L-R} is connected to the geometric properties of the intersection of a null hyperplane with a lightcone in the Minkowski spacetime. This thus in principle provides a better conceptual understanding of the result in \cite{K-L-R}.

\section{The geometry of a spacelike surface in a vacuum spacetime}\label{Geo of Surf}

We consider a spacelike surface $S$ in a vacuum spacetime $(M,g)$. On the normal bundle $\mathpzc{N}S$ of $S$, we can introduce a null frame $\{L,\underline{L}\}$. $L,\underline{L}$ are future directed null vector fields in $\mathpzc{N}S$ and satisfies
\begin{eqnarray} \nonumber
g(L,\underline{L})=-2.
\end{eqnarray}
Such choice of a null frame is not unique, since for any positive function $a$ on $S$, $\{aL,a^{-1}\underline{L}\}$ is another null frame. Associated with a null frame $\{L,\underline{L}\}$, we can define the following 2-covariant tensor fields $\chi,\underline{\chi}$ and the1-form $\zeta$ over $S$ by
\begin{eqnarray}  \nonumber
\chi(X,Y)=g(\nabla_X L,Y), \quad
\underline{\chi} (X,Y)=g(\nabla_X \underline{L},Y), \quad
\zeta(X) = \frac{1}{2} g(\nabla_X L, \underline{L})
\end{eqnarray}
for any vector fields $X,Y$ tangential to $S$.

$\chi$ is the second fundamental form of $S$ corresponding to $L$ and $\underline{\chi}$ is the second fundamental form corresponding to $\underline{L}$. We call $\zeta$ the torsion of the null frame $\{L,\underline{L}\}$. The traces $\tr \chi,\tr\underline{\chi}$ of $\chi,\underline{\chi}$ are called the null expansions and the tracefree parts $\widehat{\chi},\widehat{\underline{\chi}}$ are called the shears.

\begin{definition}
\begin{enumerate}
\item A spacelike surface $S$ is called trapped if the two null expansions are both negative, i.e. $\tr\chi<0$ and $\tr\underline{\chi}<0$.
\item A spacelike surface $S$ is called marginally trapped if one of the null expansion is identically zero and the other is negative.
\end{enumerate}
\end{definition}
 
We choose an orthogonal frame $\{e_1,e_2\}$ on the tangent bundle $\mathpzc{T}S$ such that $\{e_1,e_2,e_3=\underline{L},e_4=L\}$ is a positive oriented frame on $(M,g)$. We denote the area element of $S$ by $\slashepsilon$. Let the Latin index $a,b,\cdots$ denote $1,2$. We define the 2-covariant $S$-tensor fields $\alpha,\underline{\alpha}$, the $S$-tangential 1-forms $\beta,\underline{\beta}$ and the functions $\rho,\sigma$ on $S$ as following:
\begin{eqnarray}  \nonumber
&\alpha_{ab}=\R_{4ab4}, &\underline{\alpha}_{ab}=\R_{3ab3}, 
\\  \nonumber
&\beta_{a}=\frac{1}{2}\R_{a443}, &\underline{\beta}_{a}=\frac{1}{2}\R_{a343},
\\  \nonumber
&\rho=\frac{1}{4}\R_{4334}, &\sigma\cdot\slashepsilon_{ab}=\frac{1}{2}\R_{ab43}.
\end{eqnarray}

Because of the symmetries of the curvature tensor and the Ricci flat condition, we have
\begin{eqnarray}  \nonumber
&\R_{abc3}=-\slashepsilon_{ab}\slashepsilon_{c}^{\phantom{c}d}\underline{\beta}_d, &\R_{abc4}=\slashepsilon_{ab}\slashepsilon_{c}^{\phantom{c}d}\beta_d
\\  \nonumber
&\R_{abcd}=\rho \cdot \slashepsilon_{ab} \slashepsilon_{cd}, &\R_{3ab4}=-\rho\cdot \slash \llap g_{ab} + \sigma\cdot \slashepsilon_{ab}.
\end{eqnarray}
Hence $\{\alpha,\underline{\alpha},\beta,\underline{\beta}, \rho,\sigma\}$ contains all the curvature components.

We have the following Gauss equation:
\begin{eqnarray}
K+\frac{1}{4}\tr\chi\,\tr\underline{\chi}-\frac{1}{2}(\widehat\chi,\widehat{\underline{\chi}})=-\rho. \label{Gauss Equation}
\end{eqnarray}
Here
\begin{equation}  \nonumber
(\widehat\chi,\widehat{\underline{\chi}})= (\slash \llap g^{-1})^{ac}\,(\slash \llap g^{-1})^{bd}\, \widehat{\chi}_{ab}\,\widehat{\underline{\chi}}_{cd}.
\end{equation}

\section{The Minkowski spacetime $( \mathbb{M}, \etaup)$ and the representation of a surface in the lightcone in $\mathbb{M}$} 

The Minkowski spacetime is the 4-dim affine space $\mathbb{M}$ endowed with the flat Lorentzian metric $\etaup$. Let $\{x_0,x_1,x_2,x_3\}$ be the rectangular coordinate system on $\mathbb{M}$ such that $\etaup=- \d x_0^2 +\d x_1^2 +\d x_2^2 +\d x_3^2$. The lightcone issuing from the origin $o$ is the subset $C_o=\{- x_0^2 + x_1^2 + x_2^2 + x_3^2 =0\}$.

Let $\{ t, r, \vartheta, \phi \}$ be the spatial polar coordinate system where $ t\in \mathbb{R}$, $r\in \mathbb{R}_{>0}$, $ \vartheta \in (0,\pi)$ and $ \phi\in (0,2\pi)$. The coordinate transformation from the spatial polar coordinates to the rectangular coordinates is
\begin{equation}  \nonumber
(x_0,x_1,x_2,x_3) = (t, r \sin \vartheta \cos \phi, r \sin\vartheta \sin\phi, r \cos \vartheta).
\end{equation}

We define two optical functions $u =\frac{1}{2}(t -r) $ and $v= \frac{1}{2} (t+r)$. Then denote $\{C_u\}$ the level sets of $u$ and $\{\underline{C}_v\}$  the level sets of $v$. $C_u$ is the future half lightcone issuing from the point $(u,0,0,0)$ and $\underline{C}_v$ is the past half lightcone issuing from the point $(v,0,0,0)$. We denote $S_{u,v}$ the intersection $C_u$ with $\underline{C}_v$. 

$\{u,v,\vartheta,\phi\}$ is also a coordinate system on $\mathbb{M}$ and we call it a double null coordinate system. The metric $\etaup$ takes the form $\etaup= -2(\d u \otimes \d v + \d v \otimes \d u) + r^2 ( \d \vartheta^2 + \sin^2 \vartheta \d \phi^2)$ with $r=v-u$. The intrinsic geometry of $S_{u,v}$ is a round sphere of radius $v-u$.

The coordinate vector fields $\partial_v$ and $\partial_u$ are two geodesic null vector fields generating the lightcones $C_u$ and $\underline{C}_v$. Moreover they are normal to the surface $S_{u,v}$ and $\etaup({\partial_u},{\partial_v})=-2$. Hence $\{{\partial_u},{\partial_v} \}$ is a null frame on $S_{u,v}$. Let us denote $\{{\partial_v},{\partial_u} \}$ by $\{L,\underline{L}\}$. Then we can apply the constructions in the section \ref{Geo of Surf}. We have
\begin{equation}  \nonumber
\tr\chi=\frac{2}{r}=\frac{2}{v-u},  \quad \tr\underline{\chi}=-\frac{2}{r}= -\frac{2}{v-u}, \quad \hat{\chi}=0, \quad \hat{\underline{\chi}}=0, \quad \zeta=0,
\end{equation}
and
\begin{equation}  \nonumber
\alpha=0, \quad \underline{\alpha}=0, \quad \beta=0, \quad \underline{\beta}=0, \quad \rho=0, \quad \sigma =0.
\end{equation}

We consider the lightcone $\underline{C}_0$ which is the past half lightcone issuing from the origin. $\{u,\vartheta,\phi\}$ is a coordinate system on $\underline{C}_0$. The metric $\etaup$ restricted to $\underline{C}_0$ takes the form $\etaup|_{\underline{C}_0} = u^2(\d \vartheta^2 + \sin^2 \vartheta \d \phi^2)$.

Let $S$ be a spacelike surface in $\underline{C}_0$, we could represent it as a graph of a positive function $f$ on the sphere: $S=\{(u=-f(\vartheta,\phi),\vartheta,\phi)\}$. Conversely, for any positive function $f$ on the sphere, we could construct a spacelike surface $S_f$ in $\underline{C}_0$ by $S_f=\{(u=-f(\vartheta,\phi),\vartheta,\phi)\}$. The metric $\etaup$ restricted to $S_f$ takes the form $\etaup|_{S_f} =f^2(\vartheta,\phi)(\d \vartheta^2 + \sin^2 \vartheta \d \phi^2)$, which is a conformal deformation of the round sphere.

\section{The sections of a lightcone in the Minkowski spacetime $(\mathbb{M},\etaup)$} \label{section 4}

We investigate the geometry of the surface $S_f$ constructed above. The tangent space of $S_f$ is spanned by $\{\widetilde{\partial_\vartheta},\widetilde{\partial_\phi} \}$ where
\begin{equation}  \nonumber
\widetilde{\partial_{\vartheta}} = \partial_{\vartheta} - f_{\vartheta} \partial_u, \quad \widetilde{\partial_{\phi}}=\partial_{\phi} - f_{\phi} \partial_u.
\end{equation}
$\underline{L}={\partial_u}$ is a null vector field normal to $S_f$. By solving the equation
\begin{equation}  \nonumber
\etaup( \widetilde{L},\widetilde{{\partial_\vartheta}}) =\etaup( \widetilde{L},\widetilde{{\partial_\phi}})=0, \quad \etaup( \widetilde{L},\underline{L})= -2,
\end{equation}
we get $\widetilde{L}=\partial_v - 2\etaup^{\vartheta\vartheta}f_{\vartheta}\partial_{\vartheta} - 2\etaup^{\phi\phi} f_{\phi} \partial_{\phi}$. Hence $\{\widetilde{L},\widetilde{\underline{L}}\}$ is a null frame on the normal bundle of $S_f$ where
\begin{equation}  \nonumber
\widetilde{\underline{L}}=\underline{L}=\partial_u, \quad \widetilde{L} =\partial_v - 2\etaup^{\vartheta\vartheta}f_{\vartheta}\partial_{\vartheta} - 2\etaup^{\phi\phi} f_{\phi} \partial_{\phi}.
\end{equation}

We can calculate the null expansions $\tr{\widetilde{\chi}}$ and $\tr\widetilde{\underline{\chi}}$, the shears $\widehat{\widetilde{\chi}}$ and $\widehat{\widetilde{\underline{\chi}}}$, the torsion $\widetilde{\zeta}$ and the curvature components $\widetilde{\alpha}, \widetilde{\underline{\alpha}}, \widetilde{\beta},\widetilde{\underline{\beta}},\widetilde{\rho},\widetilde{\sigma}$ along $S_f$ corresponding to $\{\widetilde{L},\widetilde{\underline{L}}\}$. Since the Minkowski spacetime is flat, the curvature components are all zero. For our purpose, we are more concerning with the null expansions,
\begin{eqnarray}  \nonumber
& \tr \widetilde{\underline{\chi}} &=-\frac{2}{f},
\\
\label{transformation formula in Minkowski lightcone}
& \tr\widetilde{\chi} &= \frac{2}{f}  - 2 \Delta_{S_f} f + \frac{2}{f} | \d f |^2_{\etaup_{S_f}} 
\\
& &=\frac{2}{f}  (1 -  \circDelta \log f),  \nonumber
\end{eqnarray}
where $\circDelta$ is the Laplacian operator on the round sphere of radius $1$. One can easily get the formula for $\tr \widetilde{\underline{\chi}}$ from the Gauss equation \eqref{Gauss Equation}. Since the intrinsic geometry of $S_f$ is conformal to the round sphere with the conformal factor $f^2$, we have
\begin{equation}
K_{S_f} = f^{-2}( 1 - \circDelta \log f) = - \frac{1}{4} \tr\widetilde{\chi}\; \tr \widetilde{\underline{\chi}},
\end{equation}
which is consistent with the our formulas for $\tr\widetilde{\chi}$ and $\tr\widetilde{\underline{\chi}}$. 

One can not find a closed trapped spacelike section in $\underline{C}_0$, since the Minkowski spacetime is complete. However as we already mentioned in the introduction that \cite{K-L-R} proved that one could find a closed trapped surface in a vacuum spacetime which has sufficiently strong gravitation waves concentrated near some direction. It is interesting to see how the piece of the surface contained in the flat portion of the spacetime could be trapped.

It is well know that via the stereographic projection, one can conformally deform the round sphere minus a point to the flat plane. The conformal factor of the stereographic projection is $e^{-w}$ where $w$ is the Green's function for the Laplacian on the round sphere. In the spherical coordinates,
\begin{equation}  \nonumber
w=2 \log \sin \frac{\vartheta}{2}, \quad e^{-w} =\frac{1}{\sin^2 \frac{\vartheta}{2}} = \frac{2}{1-\cos \vartheta}.
\end{equation}
Since the Laplacian $\circDelta=\frac{1}{\sin \vartheta}\partial_{\vartheta} ( \sin \vartheta \partial_{\vartheta} ) + \frac{1}{\sin^2 \vartheta} \partial_{\phi}^2$, 
we have $\circDelta w = -1$ outside the point $\vartheta=0$. In the sense of distribution, $w$ satisfies
\begin{equation}  \nonumber
\circDelta w + 1 = 4\pi \delta,
\end{equation}
where $\delta$ is the Dirac function at the point $\vartheta=0$.

Considering the surface $\widetilde{S}=S_{e^{-w}}$, then its intrinsic metric $\etaup|_{\widetilde{S}}=e^{-2w}(\d \vartheta^2 + \sin^2 \vartheta \d \phi^2)$ has zero Gauss curvature. We have that
\begin{eqnarray}  \nonumber
&&\tr\widetilde{\underline{\chi}} = -2e^{w} = \cos \vartheta -1,
\\  \nonumber
&&\tr\widetilde{\chi} = 0.
\end{eqnarray}
Hence $\widetilde{S}$ is marginally trapped. Easy to see that $\widetilde{S}$ is not compact, since it goes to infinity in the lightcone at the point $\vartheta=0$.

Parameterized by $\{\vartheta, \phi\}$, we have the following map to represent $\widetilde{S}$ in the rectangular coordinates:
\begin{equation}  \nonumber
(\vartheta,\phi) \mapsto (x_0,x_1,x_2,x_3)=(\frac{-2}{1-\cos\vartheta}, \frac{2\sin\vartheta \cos\phi}{1-\cos\vartheta}, \frac{2\sin\vartheta \sin\phi}{1-\cos\vartheta}, \frac{2\cos\vartheta}{1-\cos\vartheta}).
\end{equation}
Obviously, $\widetilde{S}$ lies in the null hyperplane $H_{n}=\{x_0 + x_3 =-2\}$.  Actually $\widetilde{S}=\underline{C}_0\cap H_n$, because by solving the following equation
\begin{equation}  \nonumber
\left\{
\begin{aligned}
&x_0 +x_3 =-2,
\\
&-x_0^2 +x_1^2 +x_2^2 +x_3^2 =0,
\end{aligned}
\right.
\end{equation}
we get exactly $\widetilde{S}$.

The intersection of $\underline{C}_0$ with the spacelike hyperplane $H_s=\{x_0 = -1\}$ is a spacelike surface $S^{\prime}=\{(-1,\sin\vartheta \cos \phi, \sin\vartheta \sin\phi, \cos\phi)\}$. $S^{\prime}$ with its intrinsic metric is a round sphere of radius 1. The intersection of $\underline{C}_0$ with the timelike hyperplane $H_t = \{x_3=-1\}$ is a spacelike surface $S^{\prime\prime}= \{x_3 =-1, -x_0^2 +x_1^2 +x_2^2 =-1\}=\{ \frac{1}{\cos\vartheta}, -\tan \vartheta \cos \phi, - \tan \vartheta \sin \phi, -1)\}$. $S^{\prime\prime}$ with its intrinsic metric is a hyperbolic surface of constant Gauss curvature $-1$. $S^{\prime\prime}$ is in fact trapped, however it is not compact.

In fact the intersection of $\underline{C}_0$ with any spacelike hyperplane is intrinsically a round sphere of constant positive Gauss curvature. The intersection with any timelike hyperplane is of constant negative Gauss curvature.

\section{The short pulse initial data and the trapped surface formation theorem}

In \cite{C}, Christodoulou introduced the short pulse initial data and proved the longtime existence theorem for such initial data. If an additional lower bound of the initial data is assumed, then a trapped surface will be formed in the domain of existence of the solution. We present here the short pulse initial data and outline the argument of the trapped surface formation theorem.

Consider the characteristic initial data problem on an outgoing lightcone $C$. Let $o$ be the vertex of the lightcone and $T$ be a timelike vector at $o$. Then we can construct an affine parameter $s$ on $C$ such that the generator $L$ of $s$ satisfies that $g(L,T)=-1$ at $o$. We choose a constant $r_0>1$ and define the new affine parameter $\underline{u}=s-r_0$. We set the initial data before $\underline{u}=0$ to be trivial.  The initial data is nontrivial in the part $\{ 0\leq\underline{u}\leq \delta\}$ for some small $\delta$. Let $\psi_0$ be a smooth 2-dimensional symmetric trace-free matrix-valued function on $\mathbb{S}^2$, depending on $s\in[0,1]$, which extends smoothly by $0$ to $s\leq 0$. We set
\begin{equation}\label{short pulse ansatz}
\psi(\underline{u},\theta) = \frac{\delta^{\frac{1}{2}}}{r_0}\psi_0 (\frac{\underline{u}}{\delta},\theta), \quad (\underline{u}, \theta) \in [0,\delta] \times \mathbb{S}^2.
\end{equation}
Then the initial data in $ 0 \leq \underline{u} \leq \delta $ is given by the 2-covariant symmetric positive definite tensor desity $m$ on $\mathbb{S}^2$ depending on $\underline{u}$
\begin{equation}  \nonumber
m=\exp \psi.
\end{equation}
The incoming energy per unit solid angle in each direction in the interval $[0,\delta]$ is the integral
\begin{equation}  \nonumber
\frac{r_0^2}{8\pi} \int_0^{\delta} e \d\underline{u}
\end{equation}
where
\begin{equation}  \nonumber
e=\frac{1}{2} |\widehat{\chi} |^2 = \frac{1}{8} |\partial_{\underline{u}} m |_{m}^2.
\end{equation}
The factor $\delta^{\frac{1}{2}}$ in the ansatz \eqref{short pulse ansatz} appears because of the desire to form trapped surfaces from the focusing of incoming waves.

Now we turn to the argument of the trapped surface formation theorem in \cite{C}. By the theorem of domain of dependence, we know that the future $M_0$ of $\{\underline{u}\leq0\}$ is flat. The future boundary of $M_0$ is the incoming lightcone $\underline{C}_{\underline{u}=0}$. Let $\Gamma$ be the geodesic tangent to $T$ and $t$ be the affine parameter of $\Gamma$ generated by $T$. We define a new function $u$ along $\Gamma$ by $u=\frac{t}{2}-r_0$ and $u$ is extended as an optical function in the development which is constant along the future light cones with vertices on $\Gamma$. Set $u_0 = -r_0$. We are interested in the development of $\underline{C}_{\underline{u}=0, u_0 \leq u \leq 0}$ and $C_{u=u_0, 0 \leq \underline{u} \leq \delta}$. In the development $\underline{u}$ is extended to be an optical function which is constant along the incoming null hypersurfaces corresponding to the level sets of $\underline{u}$ on $C_{u_0}$. Then the level sets of $u$ and $\underline{u}$ give rise to a double null foliation in the development.

Let $S_{u,\underline{u}}$ be the intersection of $C_u$ with $\underline{C}_{\underline{u}}$ and $H_{-1}$ be the spacelike hypersurface foliated by $\{ S_{u,\underline{u}} \}_{u+\underline{u}=-1}$. The existence theorem guarantees that, for fixed $r_0$ and $\psi_0$, if we choose $\delta$ suitably small depending on $r$ and $\psi_0$, we can solve the characteristic initial data problem and construct a regular spacetime at least up to $H_{-1}$. We denote by $M_{-1}$ the spacetime bounded by $C_{u_0}$, $\underline{C}_0$, $\underline{C}_{\delta}$ and $H_{-1}$. In $M_{-1}$, the outmost spacelike section of the double null foliations is $S_{-1-\delta, \delta}$.

Let $k$ be a constant greater than $1$. If we assume that the incoming energy per unit solid angle in each direction in the advanced time interval $[0,\delta]$ is not less than $\frac{k}{8\pi}$, which is the uniform on $\mathbb{S}^2$ lower bound mentioned in the introduction, then for suitably small $\delta$, $S_{-1-\delta, \delta}$ is actually trapped. For $\tr\underline{\chi}$, we have on $\underline{C}_{\delta}$
\begin{equation}  \nonumber
\left\vert \tr\underline{\chi} + \frac{2}{|u|} \right\vert \leq O(\delta |u|^{-2} ). 
\end{equation}
For $\tr\chi$, we can estimate from the second variation formula along the generators of $C_{u}$ and the incoming energy along $C_u$:
\begin{eqnarray} 
D \tr\chi &=& -\frac{1}{2} (\tr\chi)^2 - |\hat \chi |^2, \label{second variation formula}
\end{eqnarray}
We integrate the second variation formula \eqref{second variation formula} along the generators of $C_{u}$,
\begin{equation} 
\tr\chi
\leq
\frac{2}{|u|} - \int_0^{\delta} |\widehat{\chi}|^2 \d \underline{u}
=
\frac{2}{|u|} - \frac{1}{|u|^2} \int_0^{\delta} |u|^2|\widehat{\chi}|^2 \d \underline{u}.
\label{integrate second variation formula}
\end{equation}
We have the following estimation
\begin{equation} 
|u|^2 |\widehat{\chi}|^2(u,\underline{u}, \theta) \geq 2|u_0|^2 e(\underline{u},\theta) - O(\delta^{-\frac{1}{2} }). \label{incoming energy in M_{-1}  prime}
\end{equation}
Then from \eqref{integrate second variation formula} and \eqref{incoming energy in M_{-1} prime}, we get the estimation of $\tr\chi$ on $\underline{C}_{\delta}$
\begin{equation}  \nonumber
\tr\chi(u,\delta,\theta) \leq \frac{2}{|u|} - \frac{2 |u_0|^2 }{|u|^2} \int_0^{\delta} e(\underline{u},\theta) \d \underline{u} + O(\delta^{\frac{1}{2}} |u|^{-2} ).
\end{equation}
Hence if the incoming energy per unit solid angle in each direction in the advanced time interval $[0,\delta]$ is not less than $\frac{k}{8\pi}$, i.e.
\begin{equation}  \nonumber
\frac{|u_0|^2}{8\pi} \int_0^{\delta} e \d\underline{u} \geq \frac{k}{8\pi},
\end{equation}
we get
\begin{equation} \label{estimate tr chi}
\tr\chi (u,\delta ,\theta) \leq \frac{2}{|u|} - \frac{2k}{|u|^{2}} + O(\delta^{\frac{1}{2}} |u|^{-2}).
\end{equation}
Then we conclude that if $\delta$ is suitably small,
\begin{equation}  \nonumber
\tr\chi(-1-\delta, \delta,\theta) \leq \frac{2+2\delta-2k}{(1+\delta)^2} + O(\delta^{\frac{1}{2}}) <0.
\end{equation}
Hence $S_{-1-\delta,\delta}$ is trapped.

\section{The connection between sections of a lightcone in $(\mathbb{M},\etaup)$ and the anisotropic criterion for formation of trapped surfaces}

In \cite{K-L-R}, it is showed that one can find a closed trapped surface $\widetilde{S}$ in $\underline{C}_{\delta}$ provided that the incoming energy per unit solid angle is large in some neighbourhood of some direction. The keys of the argument in \cite{K-L-R} are a transformation formula for the null expansion $\tr\widetilde{\chi}$ of $\widetilde{S}$, which leads to an elliptic inequality, and the solution of this elliptic inequality. We will see how the solution of the elliptic inequality in \cite{K-L-R} is connected to the geometric properties of sections of a lightcone in $(\mathbb{M},\etaup)$ derived in section \ref{section 4}. 

We present here the transformation formula and the derivation of the elliptic inequality.  $\{u,\underline{u},\theta\}$ is the double null coordinate system in $M_{-1}$. Any spacelike section $\widetilde{S}$ in $\underline{C}_{\delta}$ can be represented as the graph of a function $f$ on $\mathbb{S}^2$, i.e. $\widetilde{S} = \{ (u= -f(\theta),\delta,\theta) \}$. Then by the transformation formula for $\tr\widetilde{\chi}$ from \cite{K-L-R},
\begin{equation} \label{transformation formula}  \nonumber
\tr\widetilde{\chi} = \tr\chi - 2\Omega \Delta f - 4\Omega\eta \cdot \nabla f - 4 \Omega^2 \widehat{\underline{\chi}}(\nabla f, \nabla f) - \Omega^2 \tr\underline{\chi} | \nabla f |^2 - 8 \Omega^2 \underline{\omega} | \nabla f |^2.
\end{equation}
For a fixed function $f$ on $\mathbb{S}^2$, as $\delta \rightarrow 0$, the transformation formula becomes
\begin{equation} \label{transformation formula in the limit} 
\tr\widetilde{\chi} = \tr\chi -  \frac{2 \circDelta \log f}{f}.
\end{equation}
Let $k(\theta)$ be the incoming energy per solid angle in the direction $\theta$, i.e.
\begin{equation}   \nonumber
k(\theta)= \frac{r^2}{8\pi} \int_0^{\delta} e(\underline{u},\theta) \d\underline{u}.
\end{equation}
Then we substitute the estimation \eqref{estimate tr chi} of $\tr\chi$ into the transformation formula \eqref{transformation formula in the limit} to get
\begin{equation} \label{estimate tr tilde chi}  \nonumber
\tr \widetilde{\chi}(\theta) \leq \frac{2}{f}(1-\circDelta \log f) - \frac{2 k(\theta)}{f^2} 
\end{equation}
as $\delta \rightarrow 0$. The goal is to find $f$ such that the right hand side is less than $0$, i.e. to find the solution of the following elliptic inequality
\begin{equation}\label{elliptic inequality}
\frac{2}{f}(1-\circDelta \log f) -\frac{2 k(\theta)}{f^2} < 0.
\end{equation}

We shall use the knowledge about the sections of a lightcone in the Minkowski spacetime from section \ref{section 4} to construct a solution$f$ of \eqref{elliptic inequality}, provided that $k$ is large in some neighbourhood of some direction.

Recall the formula \eqref{transformation formula in Minkowski lightcone} in section \ref{section 4}, which is a special case of the transformation formula \eqref{transformation formula in the limit},  
\begin{equation}  \nonumber
\tr \tilde{\chi} = \frac{2}{f}(1-\circDelta \log f).
\end{equation}
This coincides with the first term of the right hand side of the elliptic inequality \eqref{elliptic inequality}. Intuitively, we can interpret the right hand side of $\eqref{elliptic inequality}$ as following: the first term $\frac{2}{f}(1-\circDelta \log f)$ comes from the intrinsic metric of $\tilde{S}$, which is like in the flat spacetime and the second term $-\frac{2 k(\theta)}{f^2} $comes from the incoming energy of the gravitational wave. The second term is a good term for the elliptic inequality, which shows that the incoming gravitational wave helps the surface becoming trapped.

From section \ref{section 4}, we know that the intersection of a null hyperplane with a lightcone is trapped. Such a trapped surface is not closed. It covers all the directions of the lightcone except one direction. We can cut this surface near the direction which is not covered and glue a piece of spacelike surface covering the direction onto the surface to make a closed surface. For such a closed surface, it will be trapped in the region coinciding with the original trapped surface. For the region near the direction not covered by the original trapped surface, with sufficient strong gravitational wave, we can make it being trapped. That is the geometric picture of the following construction of a solution of the elliptic inequality \eqref{elliptic inequality}. The construction here is almost identical to that of \cite{K-L-R}.

Let $\{\vartheta,\psi\}$ be the coordinate system on $\mathbb{S}^2$. We choose a positive smooth cutoff function $g$ supported in $[0,2]$ such that $g\equiv 1$ on [0,1] and decreases to 0 on $[1,2]$. Let $\varepsilon$ be an arbitrary small positive number. Denote $g_{\varepsilon} (\vartheta,\psi)= g(\frac{\vartheta}{\varepsilon})$. Then we construct $f_{\varepsilon}$ by
\begin{equation}  \nonumber
\log f_{\varepsilon}=
\left\{
\begin{aligned}
&-(1+\varepsilon) w, \quad \vartheta \in [2\varepsilon, \pi],
\\
&-(1-g_{\varepsilon})(1+\varepsilon)w- g_{\varepsilon}(1+\varepsilon) w_\varepsilon, \quad \vartheta\in[0, 2\varepsilon),
\end{aligned}
\right.
\end{equation}
where $w(\vartheta,\psi)=2 \log \sin \frac{\vartheta}{2}$ is the Green's function of the Laplacian on the round sphere and $w_\varepsilon =2 \log \sin \frac{\varepsilon}{2}$.

Let $k_\varepsilon$ in the following be:
\begin{equation}
k_{\varepsilon}=\sup_{\vartheta\in[0,2\varepsilon]} f_{\varepsilon} ( 1- \circDelta \log f_{\varepsilon})
\end{equation}
We choose initial data such that the incoming energy per unit solid angle $k(\vartheta,\psi)$ in $\{\vartheta\in[0,2\varepsilon]\}$ is larger than $k_\varepsilon$. In the following we prove that for $\delta$ suitable small, the surface $S_{f_\varepsilon}=\{(-f_{\varepsilon}(\theta), \delta, \theta)\}$ is trapped.

For $\vartheta \in [0,2\varepsilon]$,
\begin{eqnarray}  \nonumber
\tr \widetilde{\chi}
&\leq&
\frac{2}{f_\varepsilon}(1-\circDelta \log f_\varepsilon) - \frac{2 k(\vartheta,\psi)}{f_{\varepsilon}^2}
\\
&\leq&
\frac{2k_\varepsilon}{f_\varepsilon^2} - \frac{2 k(\vartheta,\psi)}{f_{\varepsilon}^2} \nonumber
\\
&<&
0. \nonumber
\end{eqnarray}
For $\vartheta\in [2\varepsilon,\pi]$, $\circDelta \log f_{\varepsilon}=-(1+\varepsilon)\circDelta w =1+\varepsilon$,
\begin{eqnarray}  \nonumber
\tr \widetilde{\chi}
&\leq&
\frac{2}{f_\varepsilon}(1-\circDelta \log f_\varepsilon) - \frac{2 k(\vartheta,\psi)}{f_{\varepsilon}^2}
\\
&\leq&
-\frac{2\varepsilon}{f_{\varepsilon}}  \nonumber
\\
&<&0. \nonumber
\end{eqnarray}
Hence $\tr\widetilde{\chi}$ is negative on $S_{f_{\varepsilon}}$ for $\delta$ suitable small.

For the null expansion $\tr\widetilde{\underline{\chi}}$ on $S_{f_\varepsilon}$, since $\tr\widetilde{\underline{\chi}}=\tr\underline{\chi}$
on $\underline{C}_{\delta}$,
\begin{equation}  \nonumber
\left\vert \tr\underline{\chi} + \frac{2}{|u|} \right\vert \leq O(\delta |u|^{-2} ),
\end{equation}
we get as $\delta \rightarrow 0$,
\begin{equation}  \nonumber
\tr\widetilde{\underline{\chi}} = -\frac{2}{f_{\varepsilon}}<0.
\end{equation}
Hence $\tr\widetilde{\underline{\chi}}$ is negative on $S_{f_{\varepsilon}}$ for $\delta$ suitable small. Hence $S_{f_{\varepsilon}}$ is trapped.

Now we can estimate the size of the trapped surface $S_{f_\varepsilon}$. In the direction $\vartheta=\pi$, $S_{f_{\varepsilon}}$ is near to the section $S_{-1-\delta,\delta}$ of the double null foliation, since $f_{\vartheta}=1$ at $\vartheta=\pi$. In the direction $\vartheta=0$, $S_{f_{\varepsilon}}$ is far away from the section $S_{-1-\delta,\delta}$ of the double null foliation, since $f_{\varepsilon}=\left(\frac{2}{1 -\cos \varepsilon} \right)^{1+\varepsilon}$ at $\vartheta=0$ and $\left(\frac{2}{1 -\cos \varepsilon} \right)^{1+\varepsilon}$ is much greater than 1 for small $\varepsilon$. As $\varepsilon \rightarrow 0$, $\varepsilon^2f_{\varepsilon}$ tends to a finite limit, i.e. $f_{\varepsilon} \sim \varepsilon^{-2}$.

We can also obtain a upper bound of $k_{\varepsilon}$. In the coordinate system $\{\vartheta,\psi\}$, $\circDelta= \partial_{\vartheta}^2 +\cot \vartheta \partial_\vartheta + \sin^{-2}\vartheta \partial_\psi^2$. Then for $\vartheta\in[2\varepsilon,\varepsilon]$,
\begin{eqnarray}\nonumber
w = O( \log \varepsilon), \quad 
\partial_{\vartheta} g_{\varepsilon} =O( \varepsilon^{-1}), \quad
\partial_{\vartheta} w = O( \varepsilon^{-1}),   \quad
\circDelta g_{\varepsilon} = O( \varepsilon^{-2}),  \quad
\circDelta w =1.
\end{eqnarray}
Hence $\circDelta \log f_{\varepsilon} =O( \varepsilon^{-2}\log \varepsilon)$ for $\vartheta \in [\varepsilon, 2\varepsilon]$. For $\vartheta \in [0, \varepsilon]$, $\circDelta \log f_{\varepsilon}=0$ since $\log f_{\varepsilon}$ is constant. So for $\vartheta \in [0,2\varepsilon]$, $\circDelta \log f_{\varepsilon} = O(\varepsilon^{-2}\log\varepsilon)$. Then we conclude that $k_{\varepsilon}=O(\varepsilon^{-4}\log \varepsilon)$.


\begin{thebibliography}{10}

\bibitem{C}
Christodoulou, D.,
\emph{The Formation of Black Holes in General Relativity}.
EMS Monographs in Mathematics, {\it European Mathematical Society(EMS), Z\"urich}, 2009, x+589pp,
ISBN: 978-3-03719-068-5, MR2488976, Zbl 1197.83004.

\bibitem{K-L-R}
Klainerman, S., Luk, J., Rodnianski, I.,
A fully anisotropic mechanism for formation of trapped surfaces in vacuum,
\emph{Invent. Math.} \textbf{198} (2014), no. 1, 1--26, MR3260856, Zbl 1311.35311.

\bibitem{P}
Penrose, R.,
Gravitational collapse and space-time singularities,
\emph{Phys. Rev. Lett.} \textbf{14} (1965), 57--59, MR0172678, Zbl 0125.21206




\end{thebibliography}
\end{document}